\documentclass{amsart}

\newtheorem{thm}{Theorem}

\newtheorem{prop}[thm]{Proposition}

\newtheorem{conj}[thm]{Conjecture}
   
\theoremstyle{definition}

\newtheorem{exmp}[thm]{Example}

\newtheorem{rem}[thm]{Remark}          
\newtheorem{ack}{Acknowledgments}

\newtheorem{defn-thm}[thm]{Definition--Theorem}  

\theoremstyle{remark}

\renewcommand{\c}[0]{{\mathbb C}}

\newcommand{\qtq}[1]{\quad\mbox{#1}\quad}

\newcommand{\lcm}[0]{\operatorname{lcm}}

\def\fract#1#2{\raise4pt\hbox{$ #1 \atop #2 $}}
\def\decdnar#1{\phantom{\hbox{$\scriptstyle{#1}$}}
\left\downarrow\vbox{\vskip15pt\hbox{$\scriptstyle{#1}$}}\right.}

\def\bfa{{\bf a}}

\def\bfw{{\bf w}}
\def\bfx{{\bf x}}

\def\bfz{{\bf z}}

\def\cala{{\mathcal A}}

\def\cali{{\mathcal I}}

\def\bbc{{\mathbb C}}

\def\bbp{{\mathbb P}}

\def\bbz{{\mathbb Z}}

\def\grt{\tau}

\def\grS{\Sigma}

\def\gsp1{{\mathfrak s}{\mathfrak p}(1)}

\def\la#1{\hbox to #1pc{\leftarrowfill}}
\def\ra#1{\hbox to #1pc{\rightarrowfill}}

\def\Ke{K\"ahler-Einstein }

\begin{document}
\title{Einstein Metrics on Exotic Spheres \break in Dimensions 7,11 and 15}
\author{Charles P. Boyer, Krzysztof Galicki, J\'anos Koll\'ar and Evan Thomas}
\address{CPB \& KG: Department of Mathematics \& Statistics,
University of New Mexico,
Albuquerque, NM 87131.}
\email{cboyer@math.unm.edu}
\email{galicki@math.unm.edu}
\address{JK: Department of Mathematics,
Princeton University,
Princeton, NJ 08544-1000.}
\email{kollar@math.princeton.edu}
\address{ET: Department of Physiology,
University of Melbourne,
Parkville, 3010, Australia.}
\email{evan@evan-thomas.net}

\begin{abstract}
In a recent article the first three authors proved that in dimension $4m+1$ all homotopy 
spheres that bound parallelizable manifolds
admit Einstein metrics of positive scalar curvature
which, in fact, are Sasakian-Einstein. They also conjectured that all such homotopy spheres
in dimension $4m-1, m\geq2$ admit Sasakian-Einstein metrics 
\cite{BGK}, and proved this for the simplest case, namely dimension $7.$ In this paper we 
describe computer programs that show that this conjecture is also true for 11-spheres and 
15-spheres. Moreover, a program is given that determines the partition of the 8610 
deformation classes of Sasakian-Einstein metrics into the 28 distinct oriented diffomorphism 
types in dimension $7.$
\end{abstract}

\maketitle

\section{Introduction} 
\bigskip
In a recent article the first three authors gave a 
method for constructing Einstein metrics of positive scalar curvature
on odd-dimensional homotopy spheres \cite{BGK}. By
Kervaire and Milnor \cite{ker-mil} 
and Smale \cite{smale}, for each $n\geq 5$, differentiable
homotopy spheres of dimension $n$ form an Abelian group $\Theta_n$,
where the group operation is connected sum.
$\Theta_n$ has a subgroup $bP_{n+1}$ consisting of 
those homotopy $n$-spheres which bound parallelizable manifolds $V_{n+1}.$
Kervaire and Milnor \cite{ker-mil} proved that $bP_{2m+1}=0$  for $m\geq 1,$
$bP_{4m+2}=0,$ or $\bbz_2$ and is $\bbz_2$ if $4m+2\neq 2^i-2$
for any $i\geq 3.$ The most interesting groups are
 $bP_{4m}$ for $m\geq 2.$ These are cyclic of order
$$
|bP_{4m}|=2^{2m-2}(2^{2m-1}-1)~\hbox{numerator}~\biggl(\frac{4B_m}{m}\biggr),
$$
where $B_m$ is the $m$-th Bernoulli number. 
Thus, for example $|bP_8|=28, |bP_{12}|=992,
|bP_{16}|=8128$ and $|bP_{20}|=130,816$. In the first two cases these 
include all exotic spheres.
The correspondence is given by 
$$
KM: \Sigma\mapsto \tfrac18\tau(V_{4m}(\Sigma))\!\! \mod |bP_{4m}|,
$$
where $V_{4m}(\Sigma)$ is any parallelizable manifold bounding $\Sigma$
and $\tau$ is its signature.
Let $\Sigma_i$ denote the exotic sphere with $KM(\Sigma_i)=i$.

In \cite{BGK} the authors proposed the following: 

\begin{conj} The construction of \cite{BGK} yields Einstein metrics on
every  exotic sphere that bounds a parallelizable manifold.
\end{conj}

The construction is described in sections 2--3. 
The method gives Einstein metrics whose isometry group
is one dimensional and they are even Sasakian-Einstein.

In \cite{BGK} 
the conjecture was shown to be true  in dimensions $4m+1$.
 In dimension $7$  we were also
able to verify it; the relevant signature calculations were carried out
by a computer.

The main aim of this paper is  to provide more evidence for our
conjecture by demonstrating that it is true in
dimensions 11 and 15 as well. More precisely we show 

\begin{thm}\label{1115thm} Every homotopy sphere  $\Sigma_i\in bP_{12}$ 
and $\Sigma_i\in bP_{16}$
admits at least  one 
Einstein metric.
\end{thm}

We also give a complete enumeration of all oriented diffeomorphism types in dimension $7,$ 
namely,

\begin{thm}\label{7thm} In dimension 7,  $\Sigma_i$ admits at least 
$n_i$ inequivalent deformation classes of 
Einstein  metrics,
where $(n_1,\ldots,n_{28})=(376,336,260,294,231,
284,322,402,317,309, 252,304,258,
390,409,352,\break 226,260,243,309,292,452,307,298,230,307,264,353),$ giving a total
of 
$8610$ cases.
\end{thm}

Actually for dimensions $11$ and $15,$  just as in dimension $7$, we do
get several deformation types, but the signature was computed only for a 
sample of all cases. For instance, in dimension 15 our method
gives at least $10^{50}$ deformation classes of Einstein  metrics
on all homotopy 15-spheres, and even their complete
enumeration is impossible with the current programs and facilities.

\bigskip

\section{Brieskorn-Pham Singularities and Their Links}

For $\bfa=(a_1,\ldots, a_m)\in\bbz_+^m$ 
set $F_{\bfa}(\bfz):=\sum_{i=1}^mz_i^{a_i}$. 
Consider a Brieskorn--Pham singularity
$$
Y(\bfa):=(F_{\bfa}(\bfz)=0)\subset \c^m,
\qtq{and its link} 
L(\bfa):=Y(\bfa)\cap S^{2m-1}(1).
$$
Set $C=\lcm(a_i:i=1,\dots,m)$. Both
$Y(\bfa)$ and $L(\bfa)$  are invariant under the $\c^*$-action
$$
(z_1,\dots,z_m)\mapsto (\lambda^{C/a_1}z_1,\dots,\lambda^{C/a_m}z_m).
$$
If we denote $\bfw=(w_1,\ldots,w_m)=(C/a_1,\ldots,C/a_m),$ then 
$F_{\bfa}$ is a weighted homogeneous polynomial on $\c^m$ with weight
$\bfw$ and degree $C$, i.e., 
$$F_{\bfa}(\lambda^{w_1}z_1,\dots,\lambda^{w_m}z_m)=
\lambda^CF_{\bfa}(z_1,\ldots,z_m).$$
Consider the orbit spaces: $X^{orb}(\bfa):=Y(\bfa)\setminus\{0\}/\bbc^*$
 and the
weighted projective space $\bbp(\bfw):=(\c^m\setminus\{0\}/\c^*)$.
We get a commutative diagram
\begin{equation*}
\begin{matrix}
L(\bfa) &\ra{2.5}& S^{2m-1}&\\
  \decdnar{\pi}&&\decdnar{\phantom{\pi}} &\\
  \hphantom{orb}X^{orb}(\bfa) &\ra{2.5} &\bbp(\bfw).&
\end{matrix}
\end{equation*}

It is known that the sphere $S^{2m-1}$ can be given a Sasakian structure
with respect to  the projection $S^{2m-1}\ra{1.5} \bbp(\bfw)$
  associated to
the characteristic foliation \cite{MR86g:53001}.
 In such a case the embedding $L(\bfa)\ra{1.5} S^{2m-1}$
is Sasakian and $X^{orb}(\bfa)$ is the horizontal space of the 
characteristic foliation of the link $L(\bfa)$ \cite{BG01}. 
\bigskip
\section{Orbifolds and Einstein Metrics}
\bigskip
Let $C^j=\lcm(a_i:i\neq j),\quad
b_j=\gcd(a_j,C^j)\qtq{and}  d_j=a_j/b_j.
$ 
The following result was established in 
\cite{BGK}.

\begin{thm}\label{klt}
The orbifold $X^{orb}(\bfa)=Y(\bfa)/\c^*$ is Fano and
has a \Ke metric if 
\begin{enumerate}
\item 
$1<\sum_{i=1}^m \frac1{a_i}$,
\item 
$\sum_{i=1}^m \frac1{a_i}<1+\frac{m-1}{m-2}\min_{i}\{\frac1{a_i}\}$, and
\item 
$\sum_{i=1}^m \frac1{a_i}<1+\frac{m-1}{m-2}\min_{i,j}\{\frac1{b_ib_j}\}$.
\end{enumerate}
In this case the link $L(\bfa)$ admits a Sasakian-Einstein metric
with one-dimensional isometry group.
\end{thm}

The first inequality is necessary for $X^{orb}(\bfa)$ to be Fano. 
Hence, it is also
necessary for the link $L(\bfa)$ to admit any  Sasakian-Einstein structure.
The second inequality is necessary for our algebraic approach to 
\Ke metrics to work, while the third inequality is most likely an
artifice of our estimates. Hopefully, it is not needed at all.
We should reiterate
that the failure
of our  method does not imply that $X^{orb}(\bfa)$ cannot admit
a positive K\"ahler-Einstein metric as long as $X^{orb}(\bfa)$ is Fano.

For any $m\geq3$ there are infinitely many
m-tuples satisfying the conditions of Theorem \ref{klt}. For example,
we can take $\bfa=(m-1,\ldots,m-1,k)$, where
gcd$(m-1,k)=1$ and $k>(m-1)(m-2)$. However, in this paper we are
only interested in the case when the link $L(\bfa)$ is a homotopy sphere
and, as we shall see, $L(m-1,\ldots,m-1,k)$ is not.

\bigskip
\section{Homotopy Spheres as Brieskorn-Pham Links}
\bigskip

To every m-tuple $\bfa$ one can associate a graph 
$G(\bfa)$ whose $m$ vertices are labeled by $a_1,\cdots,a_m.$ Two vertices 
$a_i$ and $a_j$ are connected if and only if $\gcd(a_i,a_j)>1.$ Let $C_{ev}$ 
denote the connected component of $G(\bfa)$ determined by the even integers. 
Note that all even vertices belong to $C_{ev},$ but $C_{ev}$ may contain odd 
vertices as well. Brieskorn shows that:

\begin{thm}\cite{briesk} \label{Brieskorngraph}
The link $L(\bfa)$ (with $m>3$) is a homotopy sphere 
if and only if either of the following 
hold:
\begin{enumerate}
\item $G(\bfa)$ contains at least two isolated points, or 
\item $G(\bfa)$ contains one odd isolated point and $C_{ev}$ has an odd number 
of vertices and for any distinct $a_i,a_j\in C_{ev},~$ $\gcd(a_i,a_j)=2.$  
\end{enumerate}
\end{thm}

We observe that, in each dimension, there are only finitely many
$m$-tuples that yield homotopy spheres and satisfy the conditions
(\ref{klt}.1--3). For that we introduce

\begin{exmp}[Euclid's or Sylvester's sequence]\label{egypt.defn}
(See \cite[Sec.4.3]{gkp} or \cite[A000058]{sequences}.)

Consider the sequence defined by the recursion relation
\begin{equation*}\label{extremseq}
c_{k+1}=c_1\cdots c_k+1=
c_k^2-c_k+1
\end{equation*}
beginning with $c_1=2.$
We call this sequence the {\it extremal sequence}. It starts as
$$
2,3,7,43,1807, 3263443, 10650056950807,\ldots,
$$
and it is easy to see (cf.\ \cite[4.17]{gkp}) that
\begin{equation*}\label{recipext}
\sum_{i=1}^m\frac1{c_i}=1-\frac1{c_{m+1}-1}= 1-\frac{1}{c_1\cdots c_m}.
\end{equation*}
In \cite{sound} it was proved that if the sum of reciprocals of $m$ natural
numbers is less than 1, then it is at most $1-1/(c_{m+1}-1)$. 
Thus, in this sense the
sequence $\{c_i\}$ is extremal.
\end{exmp}

We use the sequence $c_i$ to show that the number of  
$m$-tuples that yield homotopy spheres and satisfy the conditions
(\ref{klt}.1--3) is finite. Without loss of generality
we shall assume that the exponents are arranged in the non-decreasing
order. 

\begin{prop} \label{bounds} 
Assume that  $\bfa\in \bbz_+^m$ 
satisfies the conditions (\ref{klt}.1--2) and (\ref{Brieskorngraph}).
 Then 
$a_k< (m-k+1)(c_k-1)$, for $k=1,\ldots,m-1$ and 
$a_m<\frac{m!}{m-2}(c_m-1)$. In particular, the number  of
such  $m$-tuples is finite for each $m>3$.
\end{prop}

\begin{proof} Step 1. We first observe that
$\sum_{i=1}^{m-2}\frac1{a_i}<1$. For otherwise we would have
$$
\frac1{a_{m-1}}+\frac1{a_m}<\frac{m-1}{m-2}\cdot \frac1{a_m}<\frac2{a_m},
$$
which is impossible.

Step 2. Now, assume that $\sum_{i=1}^{k}\frac1{a_i}<1$.
Then it is also $\leq 1-1/(c_{k+1}-1)$. The remaining $m-k$ reciprocals
must sum to more than $1/(c_{k+1}-1)$, hence we obtain that
$a_{k+1}\leq (m-k)(c_{k+1}-1)$.
By Step 1 this takes care of all $a_i$ for $i\leq m-1$
and also of $a_m$ if $\sum_{i=1}^{m-1}\frac1{a_i}<1$.

Thus we are left with
Step 3: $\sum_{i=1}^{m-1}\frac1{a_i}\geq 1$. 
If equality holds there is no bound for $a_m$; however,
in this case $L(\bfa)$ is not a homotopy 
sphere, since Theorem \ref{Brieskorngraph} says that at least one of the
$a_1,\dots,a_{m-1}$  (or half of it) 
is relatively prime to the others, and this implies that we cannot
get an integer as a sum of reciprocals.
Otherwise we  have
$$
\sum_{i=1}^{m-1}\frac1{a_i} > 1+\frac1{a_1\cdots a_{m-1}}
\geq 1+\frac1{m!\cdot c_1\cdots c_{m-1}}=1+\frac1{m!\cdot (c_m-1)}.
$$
Thus we obtain that
$$
1+\frac1{m!\cdot (c_m-1)}+\frac1{a_m} < \sum_{i=1}^{m}\frac1{a_i}
\leq 1+\frac{m-1}{m-2}\cdot \frac1{a_m}.
$$
Comparing the two ends gives that
$$
a_m < \frac{m!}{m-2} (c_m-1).
$$
\end{proof}

We wrote a simple program which we call {\tt candidates.c}\footnote{The codes {\tt candidates.c}
and {\tt sig.c} as well as 
all the relevant data
files mentioned later can be downloaded
at http://www.math.unm.edu/$\tilde{\phantom{o}}$galicki/papers/codes.html.}.
This is a C code which
enumerates all  ordered $m$-tuples 
satisfying the conditions (\ref{klt}.1--3) and (\ref{Brieskorngraph}.1)
or (\ref{Brieskorngraph}.2)
in any given range
$a_i^{{\rm min}}\leq a_i\leq a_i^{{\rm max}},\ \ \ i=1,\ldots,m$, 
with the condition that
$a_1^{{\rm min}}\leq\cdots\leq a_m^{{\rm min}}$. 
In principle, for any $m\geq4$,  the program can be used
to enumerate {\bf all} $m$-tuples of this type.
However,  this is not feasible already for  $m=7$.
 On the other hand, the program has the flexibility to
``hunt" for such
$m$-tuples in any specified region of the integral lattice defined by
Proposition \ref{bounds}.

\bigskip
\section{Diffeomorphism Types --  Brieskorn, Zagier, and Hirzebruch}\label{diff}
\bigskip

By Theorem 5, 
we know when $L(\bfa)$ is a homotopy sphere. We now would like
to be able to determine  the diffeomorphism types of various links. In this
article, we are only interested in the case when $m=2k+1$.

In this case, the diffeomorphism type of a homotopy sphere 
$L(\bfa)\in 
bP_{2m-2}$ is determined \cite{ker-mil} by the signature $\grt(M)$ of a parallelizable manifold 
$M$ 
whose boundary is $\grS_\bfa^{2m-3}.$ By the Milnor Fibration Theorem \cite{mil68} we can 
take 
$M$ to be the Milnor fiber $M_\bfa^{2m-2}$ which, for links of isolated singularities coming from 
weighted homogeneous polynomials is diffeomorphic to the hypersurface  $\{\bfz\in 
\bbc^{m} ~|~F(z_1,\cdots,z_{m})=1\}.$ 

Brieskorn shows that the signature of $M^{2m-2}_\bfa$ can be written combinatorially as
\begin{equation*}\leqno{\ref{diff}.1}
\begin{split}
\grt(M^{4k}_\bfa)&= \#\bigl\{\bfx\in \bbz^{2k+1}
~|~0<x_i<a_i~\hbox{and}~0<\sum_{j=0}^{2k}\frac{x_i}{a_i}
<1~\mod 2 \bigr\}\\
& - \#\bigl\{\bfx\in \bbz^{2k+1}
~|~0<x_i<a_i~\hbox{and}~1<\sum_{j=0}^{2k}\frac{x_i}{a_i} <2~\mod 2 \bigr\},
\end{split}
\end{equation*}
where $m=2k+1$.

Using a formula of Eisenstein, Zagier (cf. \cite{hir71}) has rewritten this
formula  as: 
\begin{equation*}
\grt(M^{4k}_\bfa)=\frac{(-1)^k}{N}\sum_{j=0}^{N-1}\cot\frac{\pi(2j+1)}{ 
2N}\cot\frac{\pi(2j+1)}{2a_0}\cdots \cot\frac{\pi(2j+1)}{2a_{2k}},
\leqno{\ref{diff}.2~~~~~~~~~~~~}
\end{equation*}
where $N$ is any common multiple of the $a_i$'s.
Both formulas are quite well suited to 
computer use. We wrote a second C code which we call {\tt sig.c}. 
For any $m$-tuple with $m=2k+1=5,7,9$ {\tt sig.c} computes
the signature $\grt(\bfa):=\grt(M^{4k}_\bfa)$ and the diffeomorphism type of the link  using
either of the above formulas. Furthermore, one can use {\tt sig.c} to compute
signature and diffeomorphism type of a single $m$-tuple, or one can select an arbitrary
set of $m$-tuples ${\cali}$ and compute the signature and diffeomorphism type associated to
every $m$-tuple $\bfa\in{\cali}$. One last  feature of {\tt sig.c} is that, provided an
appropriate option is chosen, the program will start computing diffeomorphism type $g(\bfa)$ of
each $m$-tuple $\bfa\in{\cali}$ until it finds all possible oriented diffeomorphism types 
in $bP_{2m-2}$ after which it stops.
\bigskip
\section{The proofs}
\bigskip
\begin{proof} 
{\bf [Theorem \ref{7thm}]}
In dimension 7 the {\tt candidates.c} can be run in the maximal range specified by Proposition 
\ref{bounds}.
The result is exactly 8610 solutions.
These solutions become an input data file ${\cali}$ for the signature 
computation using {\tt sig.c} with
either Brieskorn or Zagier formula. In the case of 5-tuples the choice is not important. 
The signature computation takes a couple of hours on a Pentium 4 processor and the
result is a list of 8610 5-tuples $\bfa=(a_1,a_2,a_3,a_4,a_5)$ each with a number 
$g(\bfa)\in\bbz_{28}$
which determines the oriented diffeomorphism type of $L(\bfa)$. The results are contained
in the output file {\tt 7spheres.txt}.
This file can be easily sorted grouping
5-tuples with the same $g(\bfa)$ and we get the result described in Theorem \ref{7thm}.
\end{proof}

\begin{proof} {\bf [Theorem \ref{1115thm}]}
In dimension 11 the {\tt candidates.c} cannot be run in the maximal range of Proposition 
\ref{bounds}.
The complete enumeration would take too long a time. Instead, the code {\tt candidates.c} is 
used to
select   7-tuples in a specified range. This will become an input file  ${\cali}$ for the 
subsequent signature computation. One important point in selecting  ${\cali}$
is that $C={\rm lcm}(a_1,\ldots,a_7)$ 
should not be too large. The time of every individual signature computation
with {\tt sig.c} is approximately linear in $C$. Another relevant point is that 
$bP_{12}=\bbz_{992}$ so that
$|{\cali}|$ should be sufficiently large. For example, we can
ask {\tt candidates.c} to search for   7-tuples in the following range: $2\leq a_1\leq6$,  $3\leq 
a_2\leq11$
and $i+1\leq a_i\leq30$ for $i=3,4,5,6,7$. This guarantees a relatively small $C<66\cdot30^5$ for 
all
  solutions and $|{\cali}|=21,535$. One should point out that there is nothing special about the 
choice
of ${\cali}$ -- other choices can be equally successful in yielding the desired result. 
We now want to determine if we find all $g(\bfa)\in\bbz_{992}$ among $\bfa\in{\cali}$.
This is done by feeding each   7-tuple $\bfa\in{\cali}$ into {\tt sig.c} with the following 
option: the program will calculate the signature $\grt(\bfa)$ and the diffeomorphism type 
$g(\bfa)$ of each
7-tuple $\bfa\in {\cali}$ in the order specified by ${\cali}$. Any 7-tuple $\bfa\in{\cali}$
with a diffeomorphism type $g(\bfa)$
not previously found gets automatically recorded into the output file ${\cala}$. Once the program 
finds all
992 oriented diffeomorphism types it stops. The output file contains a subset of
the original input file (hopefully) containing exactly 992   7-tuples. 
All this work can done on a single PC with a Pentium 4 processor.
An example of an output file ${\cala}$ called {\tt 11spheres.txt} can be 
found at the URL mentioned in the earlier footnote.
We needed approximately 9000  7-tuples to find the 992 necessary to prove Theorem 
\ref{1115thm}.

In dimension 15 we repeat the steps outlined in the 11-dimensional case. Selecting
appropriately large data file with {\tt candidates.c} is not a problem. This can be done on 
a single PC. Given that $bP_{16}\simeq\bbz_{8128}$ one needs 
an input file $\cali$ with about 80K   9-tuples for the signature computation with {\tt sig.c}. 
More challenging problem has to do with computing signature for
these 9-tuples. To minimize computing time some care should be given to
how $\cali$ is selected. The length of a single computation varies depending
on (1) $\bfa=(a_1,\ldots,a_9)$ itself; (2) the formula used for the signature computation;
(3) the processor's speed. Also, this
is an easy parallelization task because it consists of tens of
thousand runs which are almost completely independent of each other. 
The only co-ordination that is required is to stop the process when all
$g(\bfa)\in\bbz_{8128}$ are found.

We actually generated two sets $\cali_1$, $\cali_2$ for the signature calculation. The first
set $\cali_1$ was was created by appropriately restricting the size 
of all exponents. After the calculations for $\cali_1$ were completed 
a second set $\cali_2$ was chosen to select 9-tuples with a
restricted upper bound on $C={\rm lcm}(a_1,\ldots,a_9)$.
We first used the Zagier formula \ref{diff}.2 to calculate the signature $\grt(\bfa)$ of the  
9-tuples in the selected input files $\cali_1, \cali_2$. Zagier's formula 
was chosen as this calculation is much faster 
for most individual $\bfa$'s. Exactly how much faster depends on the ratio
$a_1\cdots a_9/C$. If $C=a_1\cdots a_9$ then the Brieskorn formual \ref{diff}.1
is slightly faster. On the other hand, the problem with using the Zagier
formula for very large $C$ is that there is a large round-off error on Intel x86
processors even at maximum precision.  
When $C$ is of the order of $10^9$ this error becomes large enough that
$g(\bfa)$ is sometimes calculated incorrectly.
While this was not an issue for all   5-tuples and also for carefully selected   7-tuples
the case of 15-spheres was more of a problem. Instead of forcing the program
to do a better round-off error control with the Zagier option, we decided to do the
first calculation with the Zagier formula and then
verify all signature calculations for the candidate solution
with the Brieskorn formula \ref{diff}.1. By its nature, this formula does not have any 
round-off error. At the end we actually generated two disjoint sets of
9-tuples. One is contained in the file {\tt 15spheresA.txt}. The other one
is in {\tt 15spheresB.txt}.\end{proof}

The Zagier calculation on the first set $\cali_1$ was done at the University of Melbourne on an
IBM eServer 1350 which is a cluster of 48 2.4GHz Intel Xeon
processors.  The calculation leading to the data set {\tt 15spheresA.txt} took
approximately 9500 hours of processor time and tested nearly 70,000 candidates.
The Brieskorn verification was performed on {\tt 15spheresA.txt}
at the University of New Mexico High Performance Computer 
Center on a 256 node cluster of 733Mhz processors. This required
80,000 hours of processor time. In the case of one 9-tuple the code
calculating with the Zagier formula yielded the wrong answer:  
$g_Z(3,4,8,8,9,43,83,85,97)=3323$ while
$g_B(3,4,8,8,9,43,83,85,97)=3322$ is correct. Note that for this particular
example $C_\bfa=2,118,701,160$. It is in the $10^9$ range where
{\tt sig.c}  becomes unreliable with the Zagier option. 
An additional search for a 9-tuple with that particular oriented diffeomorphism type was 
performed so that
{\tt 15spheresA.txt} actually
contains the full set of 8128 examples. We replaced it with $\bfa=(6,6,6,6,6,10,25,59,73)$ 
which came
out of $\cali_2$. Note that here $C_\bfa=646,050$ which is smaller by 3 orders of magnitude. 

Realizing that one can do much better by a careful selection of
candidates with low $C_\bfa={\rm lcm}(a_1,\ldots,a_9)$ we 
used {\tt candidates.c} to select more ``efficient" input data set $\cali_2$. As a result
it was possible to obtain all 8128 distinct $g(\bfa)$'s calculating with the Zagier option in only
about 160 hours on the University of Melbourne facility. A 
very significant improvement indeed.
That second calculation generated
{\tt 15spheresB.txt}. The Brieskorn verification was performed on {\tt 15spheresB.txt}
at the University of Melbourne facility and it took only 1700 hours. No errors were found in the
Zagier calculation which is no
surprise: a typical $C_\bfa$ for 9-tuples of $\cali_2$ was about 3 orders of magnitude lower. 

Note that one can easily improve the ``least one" statement of Theorem \ref{1115thm} by 
repeating the
same calculation with several disjoint input files ${\cali}$. It is a simple exercise
to do it for 7-tuples and much more time consuming in the case of 9-tuples. For 9-tuples we 
actually
showed that there are at least two Sasakian-Einstein metrics on
each homotopy  sphere $\sigma^{15}_i\in bP_{16}$ as the lists {\tt 15spheresA.txt} 
and {\tt 15spheresB.txt} are disjoint. 

On the other hand, to calculate
signatures of all candidate 7-tuples and 9-tuples to get the statement similar to the one
expressed in Theorem \ref{7thm} would take thousands of years with the present technology.

\begin{rem}
It is clear that our approach  breaks down for $(2n+1)$-tuples, where $n$ is ``large enough". 
What
is exactly ``large enough" depends on several factors.
Our rough estimate indicates that assuming the same facilities and the same codes are
used it would take about 100 years to do the same calculation for 19-spheres. No doubt the
{\tt sig.c} code can be improved to calculate faster. On the other hand an average
$C$ for 11-tuples will be at least $10^3$ larger that in the 9-tuple case. In addition,
$bP_{20}=130,816$ is much bigger. Taking these two factors into account, a calculation for
19-spheres would take about $10^4$ times longer that a similar calculation for 15-spheres.
\end{rem}
\bigskip
\begin{ack}  We thank Don Zagier for discussions. We further thank
the High Performance Computer Center of the University of New Mexico and
the Advanced Research Computing Center of University of Melbourne 
for the computer resources
needed to complete this project.
CPB and KG were partially supported by the NSF under grant number
DMS-0203219 and J.K.\
was partially supported by the NSF under grant number
DMS-0200883. 
\end{ack}

\bibliographystyle{amsalpha}

\providecommand{\bysame}{\leavevmode\hbox to3em{\hrulefill}\thinspace}
\providecommand{\MR}{\relax\ifhmode\unskip\space\fi MR }
\providecommand{\MRhref}[2]{%
  \href{http://www.ams.org/mathscinet-getitem?mr=#1}{#2}
}
\providecommand{\href}[2]{#2}

\end{document}